\documentclass[11pt]{amsart}

\usepackage{amsfonts,amssymb,amsthm,amsmath}
\usepackage[english]{babel}
\usepackage[all]{xy}
\usepackage{graphics}

\newtheorem{teo}{Theorem}[section]
\newtheorem{pro}[teo]{Proposition}
\newtheorem{lemma}[teo]{Lemma}

\theoremstyle{remark}
\newtheorem{oss}{Remark}

\theoremstyle{definition}
\newtheorem{defin}{Definition}

\newcommand{\Z}{\mathbb{Z}}

\newcommand{\HHH}{\mathbb{H}}

\newcommand{\m}{\mathfrak{m}}
\newcommand{\M}{\mathfrak{M}}

\newcommand{\CCC}{\mathcal{C}}
\newcommand{\AAA}{\mathcal{A}}
\newcommand{\KKK}{\mathcal{K}}
\newcommand{\HH}{\mathcal{H}}
\newcommand{\FF}{\mathcal{F}}
\newcommand{\SSS}{\mathcal{S}}

\DeclareMathOperator{\Tor}{Tor}
\DeclareMathOperator{\lk}{lk}
\DeclareMathOperator{\st}{St}

\DeclareMathOperator{\Hom}{Hom}
\DeclareMathOperator{\projdim}{projdim}

\DeclareMathOperator{\Ht}{ht}
\DeclareMathOperator{\depth}{depth}

\newcommand{\iniett}{\hookrightarrow}

\newcommand{\ergo}{\Rightarrow}
\newcommand{\sse}{\Leftrightarrow}

\newcommand{\HR}{\widetilde{H}}  

\newcommand{\0}{\varnothing}

\begin{document}

\title[A short proof of Reisner's Theorem]{A short proof of Reisner's Theorem on Cohen-Macaulay simplicial complexes}
\author{Silvano Baggio}
\address{Dipartimento di Matematica\\
Universit\`a di Bologna\\
Piazza di Porta San Donato 5\\
40126 Bologna\\ Italy}
\email{baggio@dm.unibo.it}

\begin{abstract}We present a short proof of Reisner's Theorem, characterizing which simplicial complexes have a Cohen-Macaulay face ring.  In some cases, we can also express some homological invariants of the face ring in terms of the reduced homology of the complex.
\mbox{}

\end{abstract}

\maketitle

\section*{Introduction}
Let $\Sigma$ be a simplicial complex on the set $V=\{v_1,\dots,v_n\}$, and $k[\Sigma]$ the face ring (or Stanley-Reisner ring) over the field $k$. A theorem by Reisner (\cite[Theorem 1]{rei}, here Proposition \ref{reisner}) states that $k[\Sigma]$ is Cohen-Macaulay if and only if the reduced homology of $\Sigma$, and that of all the links of its faces, is zero, except possibly in the top degree (that is, a geometric realization $|\Sigma|$  of $\Sigma$ is a homology sphere). 

The proof presented here avoids using local cohomology, neither directly as in Reisner's paper, or indirectly, via Hochster's theorem as in \cite[5.3.9]{bh}. Here the technical difficulty relies in some homological algebra, applied to certain sheaves on $\Sigma$, with the topology induced by the (reverse) face order.    
In fact, in \cite{yuz} Reisner's Theorem is proven as a corollary of a more general theorem on the rings of sections of sheaves on posets. Our method is quite similar to his, and also to that used in \cite{bac}. Our proof is more direct and self contained, and we treat separately the case where $|\Sigma|$ is a homology manifold: in that case,  proposition \ref{main} provides a description of the homology invariants $\Tor_i^R(k[\Sigma],k)$ in terms of the reduced cohomology of $\Sigma$.   

\section{Simplicial complexes and sheaves on them}
\subsection{Notations and basic facts}
By a {\em simplicial complex} over the finite set (of {\em vertices}) $V=\{v_1\dots,v_s\}$ we mean the pair $(V,\Sigma)$, where $\Sigma$ is a set of subsets of $V$ (the {\em simplexes} or {\em faces}), such that:
\begin{align*}
 &\sigma\in\Sigma,\, \tau\subset\sigma\ergo\tau\in\Sigma; \\
 &\{v\}\in\Sigma\quad\forall\, v\in V.
\end{align*}
The {\em dimension} of a face $\sigma$ is the number of vertices of $\sigma$ minus one; so, by definition, $\dim\varnothing=-1$. The dimension of  $\Sigma$ is  $\max\{\dim\sigma\mid\sigma\in\Sigma\}$. A simplicial complex is {\em pure} if all its maximal faces have the same dimension.
\begin{defin}
 Let $R$ be a ring, and $\Sigma$ a simplicial complex with vertices $V=\{v_1,\dots,v_s\}$ and faces $\SSS\subset\mathcal{P}(V)$. The {\em Stanley-Reisner algebra}, or {\em face ring} on $R$ relative to $\Sigma$ is the $R$-algebra
$$
R[\Sigma]=\frac{R[X_1,\dots,X_s]}{I},
$$ 
where $I\subset R[X_1,\dots,X_s]$ is the ideal generated by all monomials $X_{i_1}\cdots X_{i_h}$ with $\{v_{i_1},\dots,v_{i_h}\}\not\in \SSS$.
\end{defin}
\begin{defin}(See \cite[5.3]{bh})
Let $(V,\Sigma)$ be a simplicial complex of dimension $n-1$, and let $V$ be given a total order. For each $i$-dimensional face $\sigma$ we write $\sigma=[v_0,\dots,v_i]$ if $\sigma=\{v_0,\dots,v_i\}$ and $v_0<\dots<v_i$.

The {\em augmented chain complex} of $\Sigma$ is:  
$$
\CCC(\Sigma):0\to\CCC_{n-1}\stackrel{\partial}{\to}\CCC_{n-2}\to\dots\to\CCC_0\stackrel{\partial}{\to}\CCC_{-1}\to 0
$$ 
where we set 
$$
\CCC_i=\bigoplus_{\stackrel{\sigma\in\Sigma}{\dim\sigma=i}}\Z\sigma\quad\text{and}\quad \partial\sigma=\sum_{j=0}^i(-1)^j\sigma_j
$$
for all $\sigma\in\Sigma$, and $\sigma_j=[v_0,\dots,\widehat{v_j},\dots,v_i]$ for $\sigma=[v_0,\dots,v_i]$.

The $i$-th {\em reduced simplicial homology} of $\Sigma$ with values in an abelian group $G$ is:
$$
\HR_i(\Sigma,G)=H_i(\CCC(\Sigma)\otimes G)\quad i=-1,\dots,n-1.
$$
The dual cochain complex $\Hom_{\Z}(\CCC(\Sigma),G)$ has differentials $\bar{\partial}$, defined as: $(\bar{\partial}\phi)(\alpha)=\phi(\partial \alpha)$, for  $\phi\in\Hom_{\Z}(\CCC_i,G)$, $\alpha\in\CCC_{i+1}$.   
The $i$-th group of {\em reduced simplicial cohomology} of $\Sigma$ with values in $G$ is:
$$
\HR^i(\Sigma,G)=H^i(\Hom_{\Z}(\CCC(\Sigma),G))\quad i=-1,\dots,n-1.
$$
\end{defin}
If $\sigma$ is a face of the simplicial complex $\Sigma$, the {\em link} of $\sigma$ in $\Sigma$ is $\lk_{\Sigma}\sigma=\lk\sigma=\{\tau\in\Sigma\mid \sigma\not\subset\tau,\,\sigma\cup\tau\in\Sigma\}$. It is easy to see that $\lk_{\Sigma}\sigma$ is itself a simplicial complex over the set $\{v\in V\mid v\in\tau\;\exists\,\tau\in\lk\sigma\}$.

 We denote by $\st\sigma=\{\tau\in\Sigma\mid \sigma\subset\tau\}$ the {\em star} of $\sigma$ in $\Sigma$, and by $\overline{\st}\sigma$ the least subcomplex of $\Sigma$ containing $\st\sigma$.

\begin{lemma}\label{local}Let $\Sigma$ be a simplicial complex on the vertices $\{v_1,\dots,v_t\}$, $\sigma=[v_1,\dots,v_l]$ a face of $\sigma$.   
Then we have an isomorphism of localizations of Stanley-Reisner rings:
$$
k[\Sigma]_{(X_{\sigma})}\cong k[\overline{\st}\sigma]_{(X_{\sigma})},
$$ 
where $X_{\sigma}$ is the image of the monomial $X_1\dots X_l$.
\end{lemma}
\begin{proof}
Let the vertices of  $\overline{\st}\sigma$ be $\{v_1,\dots,v_r\}$ ($l\le r\le t$). By definition $k[\Sigma]=k[X_1,\dots,X_t]/I_{\Sigma}$, where $I_{\Sigma}=(X_{i_1}\dots X_{i_k}\mid [v_{i_1},\dots,v_{i_k}]\not\in\Sigma)$, while $k[\overline{\st}\sigma]=k[X_1,\dots,X_r]/I_{\sigma}$, where $I_{\sigma}=(X_{i_1}\dots X_{i_k}\mid [v_{i_1},\dots,v_{i_k}]\not\in\overline{\st}\sigma)$. After we localize to the multiplicative system $\{X_{\sigma}^n\mid n\geq 0\}$, all monomials $X_{i_1}\dots X_{i_k}\in k[\Sigma]$ such that $[v_{i_1},\dots,v_{i_k}]\not\in \overline{\st}\sigma$ vanish: so the inclusion $k[X_1,\dots,X_r]\iniett k[X_1,\dots,X_t]$ can be lifted to a well defined ring homomorphism $ k[\Sigma]_{(X_{\sigma})}\to k[\overline{\st}\sigma]_{(X_{\sigma})}$ that is easily seen to be injective and surjective.  
\end{proof}   
\begin{lemma}\label{starlink}
With the same notation as above, 
$$
 k[\lk\sigma]=\frac{ k[\overline{\st}\sigma]}{(X_{\sigma})}.
$$ 
\end{lemma}
\begin{proof}
$
 k[\lk\sigma]= k[X_{l+1},\dots,X_r]/I_{\lk\sigma}$, with $I_{\lk\sigma}=I_{\sigma}\cap k[X_{l+1},\dots,X_r]$.
\end{proof}

\begin{pro}\label{reisner} \cite[Theorem 2]{rei}.

Let $(V,\Sigma)$ be a simplicial complex. 

The ring  $k[\Sigma]$ is  Cohen-Macaulay if and only if 
\begin{equation}\label{nullhom}
\HR_i(\lk\sigma,k)=0 \quad \forall i<\dim (\lk\sigma)\quad\forall\sigma\in\Sigma,
\end{equation}
and
\begin{equation}\label{nullhom2}
\HR_i(\Sigma,k)=0 \quad  \forall i<\dim \Sigma.
\end{equation}
\end{pro}

\begin{oss}\label{equiv}
Condition \eqref{nullhom} depends only on the topology of the support of $\Sigma$. In fact, by \cite[Prop.~4.3]{sta}, it is equivalent to: 
\[
\HR_i(|\Sigma|,|\Sigma|\smallsetminus x,k)\quad \forall i<\dim(\Sigma),\;\forall x\in|\Sigma|,
\] 
where $|\Sigma|$ is a given geometric realization of $\Sigma$.

Moreover, we can replace condition \eqref{nullhom} with 
\begin{equation}
\HR^i(\lk\sigma,k)=0 \quad \forall i<\dim (\lk\sigma)\quad\forall\sigma\in\Sigma,
\end{equation}
that is, replace (reduced) homology with cohomology. This is a consequence of the Universal Coefficient Theorems (\cite[Thm III.4.1 and Thm V.11.1]{mac}), and shall be used in the proof of Proposition \ref{main}. 
\end{oss}

\subsection{Sheaves on simplicial complexes}
Sheaves on posets have been often used to study properties of rings which can be expressed as global sections (See \cite{bac},\cite{yuz}, and also \cite{BBFK},\cite{brelu},\cite{bri}). 
  
A simplicial complex $\Sigma$ can be considered as a topological space, where the open sets are the subcomplexes. More generally, every poset $(X,\le)$ can be given a topology, where the open sets are the increasing subsets, that is, the subsets $Y\subset X$ satisfying: $y\in Y, x\in X, y\le x \ergo x\in Y$. These two topologies coincide on a simplicial complex, provided we take the the reverse face order.

Some obvious remarks: every point (face) $\sigma$ is contained in a least open subset, the subcomplex $\bar{\sigma}$ generated by $\sigma$. The closure of $\{\sigma\}$ is its star $\st\sigma=\{\tau\in\Sigma\mid \sigma\subset\tau\}$. The empty set $\varnothing$ is the maximum element in $\Sigma$ for the reverse face order, and the closure of $\{\varnothing\}$ is the whole $\Sigma$; in particular $\Sigma$ is an irreducible topological space. 

A sheaf $\FF$ of abelian groups on $\Sigma$ is by definition the data of an abelian group $\FF(\Sigma')$ ({\em sections of $\FF$ on $\Sigma'$}) for every subcomplex $\Sigma'\subset\Sigma$, and a group homomorfism ({\em restriction}) $\phi^{\Sigma_2}_{\Sigma_1}:\Sigma_2\to\Sigma_1$ for every pair $\Sigma_1\subset\Sigma_2\subset\Sigma$, such that \\
for $\Sigma_1\subset\Sigma_2\subset\Sigma_3$, we have $\phi^{\Sigma_2}_{\Sigma_1}\phi^{\Sigma_3}_{\Sigma_2}=\phi^{\Sigma_3}_{\Sigma_1}$,\\
for $\Sigma'=\cup_i\Sigma_i\subset\Sigma$, if $x_i\in\FF(\Sigma_i)$ for every $i$, and $\phi^{\Sigma_i}_{\Sigma_i\cap\Sigma_j}(x_i)=\phi^{\Sigma_j}_{\Sigma_i\cap\Sigma_j}(x_j)$ for every $i,j$, then {\em there exists a unique} $x\in\FF(\Sigma')$ such that $\phi^{\Sigma'}_{\Sigma_i}(x)=x_i$ for every $i$.    

It is clear that the stalk of $\FF$ at $\sigma\in\Sigma$ is  $\FF_{\sigma}=\FF(\bar{\sigma})$. 
The sections of $\FF$ on the subcomplex $\Sigma'$ can be described as 
\begin{equation}\label{sections}\FF(\Sigma')=\{x\in\prod_{\sigma\in\Sigma'_{max}}\FF_{\sigma}\mid {x_{\sigma}}_{|\sigma\cap\sigma'}={x_{\sigma'}}_{|\sigma\cap\sigma'}\},
\end{equation} 
where $\Sigma'_{max}$ are the maximal faces of $\Sigma'$. This implies that assigning a sheaf $\FF$ on $\Sigma$ is the same as assigning the stalks $\FF_{\sigma}$ for all $\sigma\in\Sigma$, and the restrictions $\phi^{\sigma_2}_{\sigma_1}:\FF_{\sigma_2}\to\FF_{\sigma_1}$ for every pair of faces $\sigma_1\subset\sigma_2$, with the only condition that, for $\sigma_1\subset\sigma_2\subset\sigma_3$, we have $\phi^{\sigma_3}_{\sigma_2}\phi^{\sigma_2}_{\sigma_1}=\phi^{\sigma_3}_{\sigma_1}$.  

The simplest sheaves on $\Sigma$ can be defined in the following way on the stalks. If $G$ is an abelian group, and $\sigma\in\Sigma$: $G(\sigma)_{\sigma}=G$, 
while $G(\sigma)_{\tau}=0$ if $\tau\not=\sigma$; all restrictions are zero.  We call $G(\sigma)$ the {\em simple} sheaf with support in $\sigma$ and values in $G$. The cohomology of such sheaves can be described directly in terms of the reduced cohomology of the links in $\Sigma$.
\begin{lemma}\label{semplici}
(See \cite[Lemma 3.1]{bac})
\begin{itemize}
\item[(i)] The global sections of the simple sheaf $G(\sigma)$ are
\begin{displaymath}
G(\sigma)(\Sigma)=H^0(\Sigma,G(\sigma))=\begin{cases}G\quad\text{if $\sigma$ is a maximal cone,}\\0\quad\text{otherwise.}\end{cases}
\end{displaymath}
\item[(ii)] If $i\geq 1$, then $H^i(\Delta,G(\sigma))\cong \HR^{i-1}(\lk\sigma),G)$.
\end{itemize}
\end{lemma}
\begin{proof}
(i) follows from \eqref{sections}. Let us prove (ii). The map $j_{\sigma}:\lk\sigma\to\Sigma$, sending $\tau\to\tau\cup\sigma$ is a continuous injection, and its image is $\st\sigma$, which is closed in $\Sigma$. Since $G(\sigma)$ is the push forward via $j_{\sigma}$ of the sheaf $G(\varnothing)$ on $\lk\sigma$, we have: $H^i(\Sigma,G(\sigma))=H^i(\lk\sigma,G(\varnothing))$. We need only to prove that $H^i(\Sigma,G(0))=\HR^{i-1}(\Sigma,G)$ for $i>0$.

Let $\tilde{G}$ the constant sheaf on $\Sigma$, with values in $G$. $G(\0)$ is a subsheaf of $\tilde{G}$. Let $G_{\0}=\tilde{G}/G(\0)$. 
 Since $\tilde{G}$ is acyclic, the short exact sequence
$
0\to G(\0)\to\tilde{G}\to G_{\0}\to 0
$
induces the long exact sequence in cohomology (assume $\dim\Sigma>0$):
\[
0\to H^0(\Sigma,\tilde{G})\to H^0(\Sigma,G_{\0})\to H^1(\Sigma,G(\0))\to H^1(\Sigma,\tilde{G})=0,
\]
which implies $H^1(\Sigma,G(\0))\cong \HR^0(\Sigma,G)$; moreover,
$ H^{i-1}(\Sigma,G_{\0})\cong H^i(\Sigma,G(0))$ for $i\geq 2$.
We can conclude by applying the following lemma.
\end{proof}
\begin{lemma}\cite[Theorem 2.1]{bac} 
$$H^i(\Sigma\smallsetminus \{\0\},\tilde{G})\cong\HR^i(\Sigma,G),$$ 
where $\tilde{G}$ is the constant sheaf with values in $G$, on  $\Sigma\smallsetminus \{\0\}$. 
\end{lemma}
\begin{proof} The set $\SSS=\{C_j=\st(v_j)\mid j=1,\dots,n\}$ is a closed covering of $\Sigma\smallsetminus\{\0\}$. 
For any $(p+1)$-uple of indices $i_0,\dots,i_p$, the intersection $C_{i_0}\cap \dots\cap C_{i_p}$  is either the empty set (if $[v_0\dots v_p]\not\in\Sigma$), or the star $\st[v_0\dots v_p]$. Let $G_{i_0\dots i_p}=\tilde{G}_{|C_{i_0}\cap\dots\cap C_{i_p}}$: these (constant) sheaves  are all flabby and then acyclic. 

Consider the complex
\begin{equation}\label{risoluz}
0\to \tilde{G}\to \bigoplus_{i=1}^m G_i\to \bigoplus_{1\leq i_0<i_1\leq m}G_{i_0i_1}\to \dots \to \bigoplus_{1\leq i_0<\dots<i_{n}\leq m}G_{i_0\dots i_{n}}\to 0,  
\end{equation}
where differentials are defined as follows (indices with \ $\widehat{}$\,\, are omitted):
\begin{gather}
\label{differ}
d\big((a_{i_0\dots i_k})_{i_0\dots i_k}\big) = \Big(\sum_{h=0}^{k+1} (-1)^h \bar{a}_{j_0\ldots\widehat{j_h}\dots j_{k+1}}\Big)_{j_0\dots j_{k+1}}.
\end{gather}
 The above notation means: if $a\in G_{j_0\dots \widehat{j_h}\dots j_{k+1}}(\Sigma')$, with $\Sigma'\subset\Sigma$, then $\bar{a}$ is the image of $a$ via the projection $G_{j_0\dots \widehat{j_h}\dots j_{k+1}}(\Sigma')\to G_{j_0\dots  j_{k+1}}(\Sigma')$. 

The complex \eqref{risoluz} is exact, and so it is an acyclic resolution of the sheaf $\tilde{G}$ (To see this, check that the complex of the stalks relative to each face of $\Sigma$ is exact). Since $H^0(\st(\sigma),\bar{G}_{|\st(\sigma)})=G$ for any $\sigma$, the complex of the global sections of \eqref{risoluz} is 
$$
0\to \bigoplus_{\sigma\in\Delta_1}G_{(\sigma)}\to\dots \to\bigoplus_{\sigma\in\Delta_n}G_{(\sigma)}\to 0.
$$ 
If we define  differentials as in \eqref{differ}, this is the cochain complex of $S_{\Sigma}$ with values in $G$, and  its cohomology is the reduced simplicial cohomology of $S_{\Sigma}$.
\end{proof}

\subsection{Definition of the sheaf $\AAA$}

Let $\Sigma$ be a simplicial complex of dimension $d$ over the set $V=\{v_1,\dots,v_n\}$. From now on, let $R=k[X_1,\dots,X_n]$ be the ring of the polynomials in n indeterminates on the field k.

Let us define the sheaf of $R$-algebras $\AAA$ over $\Sigma$: its sections over the subcomplex $\Sigma'$ are $\AAA(\Sigma')=R/I_{\Sigma'}$, where $I_{\Sigma'}=(X_{i_1}\dots X_{i_k}\mid [v_{i_1},\dots,v_{i_k}]\not\in\Sigma')$. If $\Sigma_1'\subset\Sigma_2'$, one can define a surjective homomorphism of $R$-algebras $R/I_{\Sigma_2'}\to R/I_{\Sigma_1'}$, sending $X_i$ to itself if $v_i\in\Sigma_1'$, to zero otherwise. This homomorphism is the restriction morphism $\AAA(\Sigma_2')\to\AAA(\Sigma_1')$. 

Equivalently, $\AAA$ can be defined as the only sheaf on $\Sigma$, such that its stalk at $\sigma=[v_{i_1},\dots,v_{i_h}]$ is $\AAA_{\sigma}=k[X_{i_1},\dots,X_{i_h}]$, and, if $\tau=[v_{i_1},\dots,v_{i_j}]$, with $j<h$ the map $\AAA_{\sigma}\to\AAA_{\tau}$ (that is,$\AAA(\sigma)\to\AAA(\tau)$ ) sends $X_{i_l}$ to itself for $l\le j$, to zero otherwise.

Obviously  $\AAA(\Sigma)=k[\Sigma]$. $\AAA$ is flasque, for the map $\AAA(\Sigma)\to\AAA(\Sigma')$ is surjective for every subcomplex $\Sigma'$.  
%
%
\section{The main proof}
First we prove Reisner's Theorem in a particular case, when the simplicial complex is a homology manifold. The proof of the general case is at the end of this section.
\begin{pro}\label{main}
Let $\Sigma$ be a pure $d$-dimensional simplicial complex on $n$ vertices, satisfying condition \eqref{nullhom}. Then 
\[
\Tor_i^R(k[\Sigma],k)=\bigoplus_{i=0}^n\HR^{i-r-1}(\Sigma,\wedge^i k^n).
\]
In particular, the reduced homology of $\Sigma$ vanishes in degree less than $d$ if and only if 
\[
\Tor_i^R(k[\Sigma],k)=0\quad \forall i>n-d-1,
\]
and this is equivalent to: $k[\Sigma]$ is a Cohen-Macaulay ring. 
\end{pro}
\begin{proof}Let $K^*(X_1,\dots,X_n)$ be the Koszul complex relative to $X_1,\dots,X_n\in R$: it is a free resolution of $k=R/(X_1,\dots,X_n)$ as an $R$-module. We can consider $K^*(X_1,\dots,X_n)$ as a complex of constant sheaves on $\Sigma$, with negative degree, from $-n$ to $0$. Let $\KKK$ be the complex of sheaves of $R$-algebras, obtained by tensoring $K^*(X_1,\dots,X_n)$ with the sheaf $\AAA$. Every $\KKK^{-i}\cong\wedge^{i} \AAA^n$ is flasque and therefore acyclic. We have: $\Tor_i^R(k[\Sigma],k)=H^{-i}(\KKK^*(\Sigma))=\HHH^{-i}(\Sigma,\KKK^*)$ for $i=1,\dots,n$, where $\HHH$ denotes the hypercohomology functor.

Let us consider the decreasing sequence of open sets $\Sigma^i=\{\sigma\in\Sigma\mid \dim\sigma\le d-i\}\stackrel{f_i}{\iniett}\Sigma$. Define $\KKK^*_i$ as the complex ${f_i}_{!}(\KKK^*_{|\Sigma^i})$, that is, $\KKK$ restricted to $\Sigma^i$ and then extended by zero to $\Sigma$. We obtain a sequence of sheaf complexes $0=\KKK^*_{d+1}\iniett\KKK^*_d\iniett\dots\iniett\KKK^*_0=\KKK^*$. Note that the quotient $\KKK^*_p/\KKK^*_{p+1}$ is supported in $\Sigma_p\setminus\Sigma_{p+1}$, which is a discrete topological space, so we can express $\KKK^*_p/\KKK^*_{p+1}$ as a direct sum of complexes of simple sheaves $\bigoplus_{\dim\sigma=d-p}\KKK^*_{\sigma}(\sigma)$. 

Standard arguments of homological algebra provide us with a spectral sequence $E_1^{pq}=\HHH^{p+q}(\Sigma,\KKK^*_p/\KKK^*_{p+1} )\ergo\HHH^{p+q}(\Sigma,\KKK^*)$. 

We can decompose the $E_1$ terms as $E_1^{pq}=\bigoplus_{\dim\sigma=d-p}\HHH^{p+q}(\Sigma,\KKK^*_{\sigma}(\sigma))$.

\emph{Claim}: condition \eqref{nullhom} implies that $E_1^{pq}=0$ for $0\le p\le d$ and $p+q<d-n+1$, or, equivalently, that, if $\dim\sigma\ge 0$, then $\HHH^l(\Sigma,\KKK^*_{\sigma}(\sigma))=0$ for $l<d-n+1$. 

\emph{Proof of the claim.} There is another standard spectral sequence, converging to the hypercohomology: $E^{ij}_2=H^i(\Sigma,\HH^j(\KKK^*_{\sigma}(\sigma)))\ergo \HHH^{i+j}(\Sigma,\KKK^*_{\sigma}(\sigma))$, where $\HH^j$ denotes the $j$-th cohomology sheaf of a complex of sheaves.  
  Lemma \ref{semplici} implies: $H^i(\Sigma,\HH^j(\KKK^*_{\sigma}(\sigma)))=\HR^{i-1}(\lk\sigma,\Tor_{-j}^R(\AAA_{\sigma},k))$. Because of condition \eqref{nullhom}, these groups vanish for $0<i<\dim(\lk\sigma)+1\sse 0<i<d-\dim\sigma$ (and they vanish trivially for $i>\dim(\lk\sigma)+1=d-\dim\sigma$). Two cases remain: $i=0$ and $i=d-\dim\sigma$.  If $\sigma$ is maximal, $\dim\sigma=d$, $\projdim_R(\AAA_{\sigma})=n-d-1$, and, by Lemma \ref{semplici}, $H^0(\Sigma,\HH^j(\KKK^*_{\sigma}(\sigma)))=H^{d-\dim\sigma}(\Sigma,\HH^j(\KKK^*_{\sigma}(\sigma)))=\Tor^R_{-j}(\AAA_{\sigma},k)=0$ for $j<d-n+1$. If $\sigma$ is not maximal, we know, by Lemma \ref{semplici}, that $H^0(\Sigma,\HH^j(\KKK^*_{\sigma}(\sigma)))=0$. Since $\Tor_{-j}^R(\AAA_{\sigma},k)=0$ for $j<-\projdim_R(\AAA_{\sigma})=\dim\sigma-n+1$, we conclude that, also when $i=d-\dim\sigma$, $E_2^{ij}=0$ for $i+j<(d-\dim\sigma)+(\dim\sigma-n+1)=d-n+1$. \emph{The claim is proved}.

So $\HHH^{r}(\Sigma,\KKK^*)=\HHH^r(\Sigma,\KKK^*_{\0})$. Since the differentials of $\KKK^*_{\0}$ are zero, we can decompose $\KKK^*_{\0}=\oplus_{i=-n}^0\KKK^*_{\0 i}$; here $\KKK^*_{\0 i}$ is the complex everywhere zero except in degree $i$, where it is the constant sheaf $\KKK^i_{\0 i}=\wedge^{-i} k^n$. Since $\HHH^r(\Sigma,\KKK^*_{\0 i})=H^{r-i}(\Sigma,\KKK^i_{\0 i})=\HR^{r-i-1}(\Sigma,\wedge^{-i}k^n)$ for $r\ge i$, we can conclude: $\HHH^r(\Sigma,\KKK^*)=\bigoplus_{i=-n}^0\HHH^r(\Sigma,\KKK^*_{\0 i})=\bigoplus_{i=-n}^0\HR^{r-i-1}(\Sigma,\wedge^{-i}k^n)=\bigoplus_{i=0}^n\HR^{i+r-1}(\Sigma,\wedge^i k^n)$ 

The last statement follows from Lemma \ref{A}.
\end{proof}
\begin{lemma}\label{A}
Let $R=k[X_1,\dots,X_n]$; let $S$ be an $R$-algebra, finitely generated as an $R$-module, and such that $\Ht\M$ is the same for all maximal ideals $\M\subset S$. Then $S$ is Cohen-Macaulay if and only if $\Tor_i^R(S,k)=0$ for all $i>n-d$, where $d=\dim S$.   
\end{lemma}
\begin{proof} By definition, $S$ is CM if and only if $d=\dim S_{\M}=\depth S_{\M}$ for all maximal ideals $\M\subset S$. The Auslander-Buchsbaum formula (\cite[Theorem 1.3.3]{bh}), applied to $S_{\M}$ as an $R_{\m}$ module (with $\m=\M\cap R$) yields: $\depth_{R_{m}}(S_{\M})+\projdim_{R_{\m}}S_{\M}=\depth R_{\m}=n$. By \cite[Theorem 1.3.2]{bh}, $\projdim_{R_{\m}}S_{\M}=\sup\{i\mid \Tor_i^{R_{\m}}(S_{\M},k)\not=0\}$. We conclude by noting that this last number equals $\sup\{i\mid \Tor_i^{R}(S,k)\not=0\}$, since $\Tor_i^{R_{\m}}(S_{\M},k)\cong R_{\m}\otimes_{R}\Tor_i^{R}(S,k)$; and $\depth_{R_{\m}}(S_{\M})=\depth_{S_{\M}}(S_{\M})$ (\cite[Ex.~1.2.26]{bh}).
\end{proof}

We give now the proof of Reisner's Theorem for a generic simplicial complex $\Sigma$.

\begin{proof}(of Proposition \ref{reisner})
 First, we recall that a Cohen-Macaulay simplicial complex is pure (\cite[Cor.~5.1.5]{bh}). So, if $k[\Sigma]$ is Cohen-Macaulay, we can apply Proposition \ref{main} and obtain conditions \eqref{nullhom} and \eqref{nullhom2}. Conversely, suppose  first that condition \eqref{nullhom} holds for $\Sigma$, but condition \eqref{nullhom2} does not. Then we can apply Proposition \ref{main}: $k[\Sigma]$ is not Cohen-Macaulay. If condition \eqref{nullhom} does not hold, then there exists a maximal  face $\sigma\in\Sigma$ among the faces such that  $\HR_i(\lk\sigma,k)\not=0$ for some $i<\dim(\lk\sigma)$. Then $\lk\sigma$ satisfies condition \eqref{nullhom} but not condition \eqref{nullhom2}. We can apply Proposition \ref{main} to $\lk\sigma$: $k[\lk\sigma]$ is not Cohen-Macaulay. By Lemmas \ref{local} and \ref{starlink}, $k[\Sigma]$ is not Cohen-Macaulay. 
\end{proof}

\end{document}